\documentclass[conference, 11pt]{IEEEtran}

\usepackage{cite}
\usepackage{amsmath,amssymb,amsfonts}
\usepackage{graphicx}
\usepackage{booktabs}
\usepackage{url}
\usepackage{hyperref}

\begin{document}

\title{A Robust Electric Vehicle Charging Station Optimization under Rare Event Disruptions}

\author{
\IEEEauthorblockN{Yan Pratama Akhra}
\IEEEauthorblockA{\textit{Industrial and Systems Engineering} \\
\textit{Institut Teknologi Sepuluh Nopember}\\
\textit{(ITS)}, Surabaya, Indonesia \\ pratamaakhra@gmail.com}
\and
\IEEEauthorblockN{Iwan Vanany}
\IEEEauthorblockA{\textit{Industrial and Systems Engineering} \\
\textit{Institut Teknologi Sepuluh Nopember}\\
\textit{(ITS)}, Surabaya, Indonesia \\ vanany@ie.its.ac.id}
\and
\IEEEauthorblockN{Mansur M. Arief}
\IEEEauthorblockA{\textit{Industrial and Systems Engineering} \\
\textit{King Fahd University of Petroleum} \\ \textit{and Minerals (KFUPM)}, Saudi Arabia \\
mansur.arief@kfupm.edu.sa}
}

\maketitle

\begin{abstract}
The rising demand for electric vehicles (EVs) requires robust and accessible charging infrastructure, particularly in developing countries where frequent electricity disruptions pose a significant service-reliability risk. Prior charging-infrastructure optimization studies rarely address such disruptions rigorously, yielding suboptimal and fragile designs. We present a simulation-based optimization model that estimates each candidate station's service reliability and embeds it directly in the objective and constraints, framing charging-network design as an engineering-management decision under uncertainty. To make the approach computationally tractable, we apply the control variates (CV) variance-reduction technique to the reliability estimation. Using a real dataset from Surabaya, Indonesia, our numerical experiment shows that the proposed robust model achieves a 13\% higher average objective value than the non-robust benchmark and remains robust even when station reliability is under- or over-estimated. The CV technique reduces the required simulation sample size by up to a factor of 10 relative to Monte Carlo, allowing the model to be solved efficiently with a standard mixed-integer programming (MIP) solver.
\end{abstract}

\begin{IEEEkeywords}
electric vehicle charging, facility location, robust optimization, simulation, variance reduction, decision under uncertainty
\end{IEEEkeywords}

\section{Introduction}

The growing global demand for electric vehicles (EVs) has brought to the forefront the need for reliable and accessible EV charging infrastructure. As governments set ambitious electrification targets, worldwide EV demand has grown rapidly, from roughly 17,000 EVs on the road in 2010 to more than one million cars sold in China alone in 2019, and the trend is projected to continue \cite{evoutlook}. Developing countries are also promoting EV adoption alongside greener electricity \cite{lai2022critical}, though consumer acceptance and adoption barriers remain decisive for uptake \cite{widiawati2022acceptance,dabic2025drone}. Indonesia, for example, targets 20\% electric automobile sales by 2025 and fully electrified transportation by 2050 \cite{indonesiaenergy}, yet had only 240 public charging points as of 2021 \cite{ev_240} against an estimated 31,000 required nationwide \cite{ev_31000}.

Designing a convenient and reliable charging network is a complex engineering-management problem. Convenient siting must account for population density and potential EV demand distribution \cite{khalid2019comprehensive}, while suitable land is scarce in dense cities. More critically, in developing countries service uncertainty---especially electricity supply---is a dominant risk. Outages and other disruptions occur frequently and undermine the service quality that smart-charging strategies \cite{fachrizal2020smart} presume. Managing this disruption risk is therefore central to any credible infrastructure plan, placing the problem squarely within operations research, decision analysis, and risk management for industrial and service systems---themes central to engineering management, from AI-enabled decision systems \cite{bleisinger2025role} and the engineering competencies they demand \cite{yu2025digital} to supply-chain disruption and resilience \cite{kiessling2025resilience,mwenya2025supply}.

\begin{figure}[t]
    \centering
    \includegraphics[width=\linewidth]{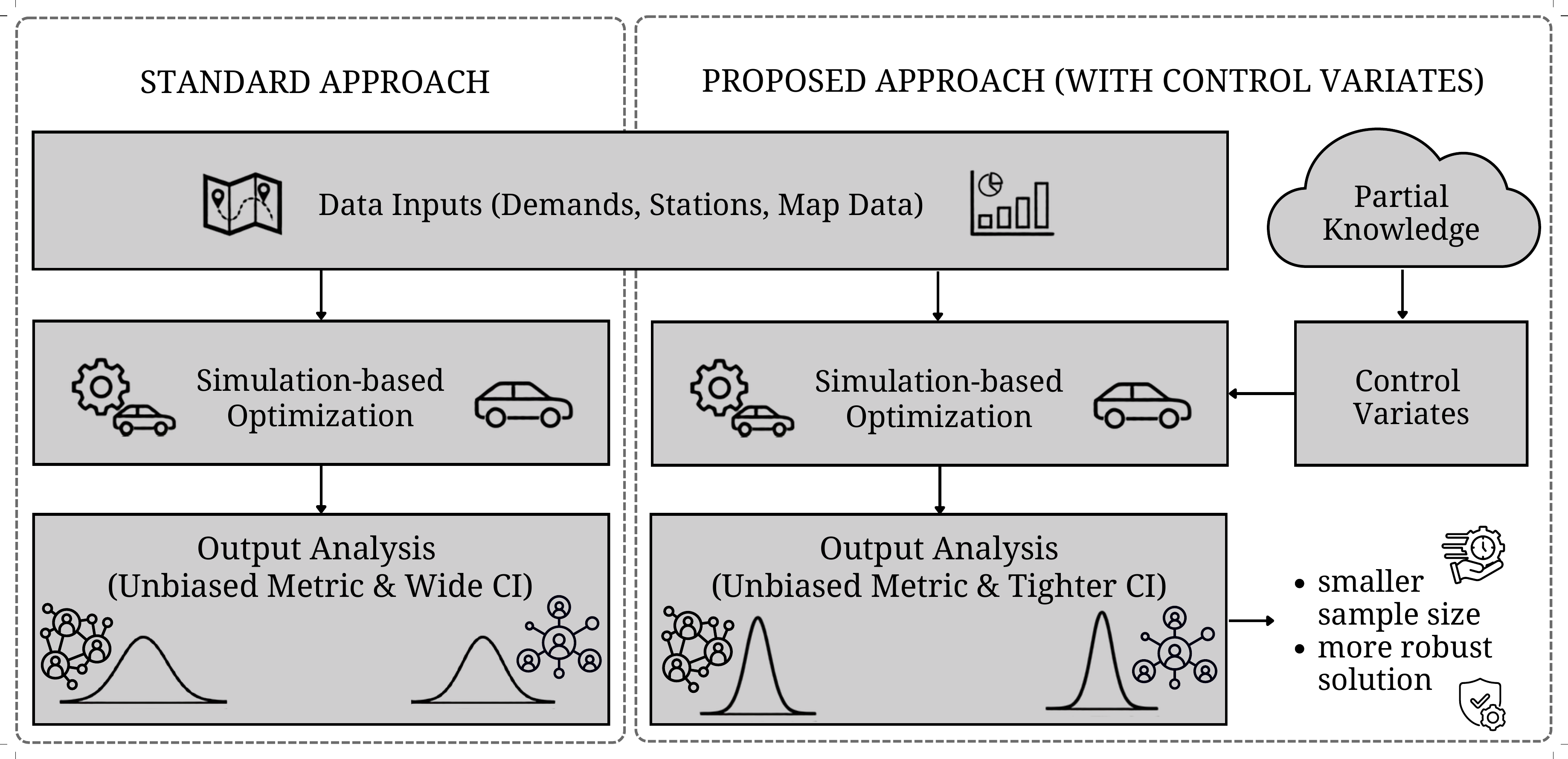}
    \caption{The proposed framework. Incorporating partial knowledge of real-world uncertainty reduces the variance of the simulation outputs, lowering the sample size needed for tighter confidence intervals.}
    \label{fig:framework}
\end{figure}

To address this challenge, we propose a robust simulation-based optimization model that estimates the service reliability of candidate charging stations and incorporates this information into both the objective function and the constraints. Unlike prior works that assume a known disruption-probability model, our simulation-based approach is versatile and data-driven. We further employ the control variates (CV) variance-reduction technique to improve simulation efficiency, reducing the required sample size by up to a factor of 10 relative to naive Monte Carlo (MC) simulation. The result is an efficient mixed-integer programming (MIP) model that balances the total cost of investing in and operating the network against the quality of service delivered to the public. Fig.~\ref{fig:framework} contrasts the traditional modeling approach without variance reduction with the proposed framework, which attains a tighter confidence interval (and thus more precise output) at a lower computational cost.

Our work contributes in three ways. First, we propose a model that directly addresses electricity-disruption risk in EV charging-station planning for developing countries. Second, we integrate the estimation of disruption probabilities into the model, giving a data-driven alternative to approaches that assume such probabilities a priori. Third, using an EV dataset from Surabaya, Indonesia, we demonstrate the robustness of the model against a non-robust benchmark, even when disruption probabilities are slightly under- or over-estimated, achieving a 13\% higher average objective value.

The rest of the paper is organized as follows. Section~\ref{sec:lit_review} reviews related work. Section~\ref{sec:problem_formulation} presents the model and the CV-based reliability estimation. Section~\ref{sec:exp} describes the experiment, Section~\ref{sec:discussion} discusses findings, and Section~\ref{sec:conclusion} concludes.

\section{Related Work}
\label{sec:lit_review}

EV charging-infrastructure planning is commonly cast as a facility location problem that minimizes cost subject to service-performance constraints. Early studies adopted deterministic models focused on minimizing the number of stations or development cost \cite{he2016sustainability} or infrastructure and access cost \cite{frade2011optimal}; related energy-replenishment siting includes $p$-median optimization of EV battery-swapping facilities \cite{li2022cabinet}. Many subsequent works continued this deterministic trend across diverse objectives, e.g., \cite{amilia2022designing}. Another line integrates charging into smart-grid or renewable-energy design, which is valuable but often impractical for dense urban areas in developing countries. A comprehensive review \cite{unterluggauer2022electric} emphasizes the need for models that account for randomness and variability, yet rigorous real-world case studies with uncertainty quantification remain scarce.

Case studies span both developed and developing regions, including slow charging in Lisbon \cite{frade2011optimal} and a city-scale simulation in Singapore \cite{bi2017simulation}. Closest to our setting, \cite{amilia2022designing} proposed a set-covering model for Surabaya but ignored electricity disruptions, relying only on redundant coverage as a buffer and yielding an oversimplified, suboptimal design. Disruption and uncertainty are addressed in only a handful of works, e.g., multi-criteria decisions against flooding \cite{zhang2023optimal}, two-stage stochastic programming with photovoltaic and storage systems \cite{liu2022optimal}, resilience-aware storage sizing \cite{hussain2020optimal}, disruption-aware location--allocation with sourcing recovery \cite{khoirani2022disruption}, and other stochastic formulations \cite{hosseini2015refueling,yildiz2019urban}. A recurring caveat is that such models become extremely hard to solve when high-confidence solutions are required.

Our work extends stochastic modeling by introducing control variates \cite{hesterberg1998control}, a variance-reduction technique, to accelerate the simulation-based optimization. Simulation lets the modeler tune fidelity to the available prior knowledge and verify the model by comparing estimated disruption probabilities against historical data. Combining stochastic modeling with control variates distinguishes this work from prior research and offers a more accurate and practical route to disruption-aware charging-infrastructure design.

\section{Model Formulation}
\label{sec:problem_formulation}

\subsection{Decision Variables}
We consider a set of demand nodes $I$ (sub-district centers) and supply nodes $J$ (candidate station locations), and $K$ vehicle types (here electric motorcycles and cars). The average travel time from node $i \in I$ to node $j \in J$ is $d_{ij}$, and $d_{max}$ upper-bounds the acceptable time-to-charge. The binary variable $x_j$ indicates whether candidate $j$ is selected, and $y_{ij}$ indicates whether demand node $i$ is assigned to station $j$. The integer variables $v_{ij}^k$ and $u_j$ denote the number of type-$k$ vehicles from node $i$ charged at node $j$ and the number of connectors installed at $j$, respectively.

Each opened station $j$ incurs a daily cost $h_j$ and can host at most $q_j$ connectors due to limited space. Each connector incurs daily operating cost $g$ and has daily throughput $c_j$ kWh. A type-$k$ vehicle draws $e_k$ kWh and takes $t_k$ minutes to fast-charge. Electricity price $r$ converts energy to revenue, and $s$ is the penalty rate for unmet demand.

\subsection{Objective Function}
The objective maximizes daily profit under random disruptions: revenue from undisrupted stations minus operating and investment costs, with a penalty on unmet demand. Each station $j$ has reliability
\begin{equation}
p_j = \mathbb{P}(Z_j \leq z_j) = \mathbb{E}\!\left[\mathbb{I}(Z_j \leq z_j)\right],
\end{equation}
where $Z = [Z_j]_{j \in J} \sim q$ models the underlying state (e.g., extreme weather, spiking demand, or fallen trees \cite{10causespoweroutages}) that triggers a disruption at station $j$ when it exceeds threshold $z_j$. When $p_j$ is known a priori, the objective is
\begin{align}
\max~&\sum_{i \in I}\sum_{j \in J}\sum_{k \in K}
\underbrace{r e_k p_j v_{ij}^k}_{\text{revenue}} - \underbrace{s d_{ij}(1-p_j)v_{ij}^k}_{\text{penalty}} \nonumber\\
&- \sum_{j \in J}\underbrace{(g u_j + h_j x_j)}_{\text{total cost}}.
\end{align}
When $p_j$ is unknown, we estimate the expectations by Monte Carlo (MC). Drawing $n$ i.i.d. samples $Z_{j1},\dots,Z_{jn}$ and computing binary disruption indicators $I_{jl} = \mathbb{I}(Z_{jl} \leq z_j)$, the robust objective becomes
\begin{align}
\max~&\sum_{i \in I}\sum_{j \in J}\sum_{k \in K}\sum_{l=1}^n \frac{1}{n}
\Big(\underbrace{r e_k v_{ij}^k I_{jl}}_{\text{revenue}} - \underbrace{s d_{ij} v_{ij}^k (1-I_{jl})}_{\text{penalty}}\Big) \nonumber\\
&- \sum_{j \in J}\underbrace{(g u_j + h_j x_j)}_{\text{total cost}}.
\label{obj_est_sim}
\end{align}
We call this the \textit{Robust Model}. Setting $I_{jl}=1$ for all $j,l$ recovers the original \textit{Non-Robust} model of \cite{amilia2022designing}. Both solutions are evaluated under disruptions generated with a different random seed.

\subsection{Constraints}
The maximization of \eqref{obj_est_sim} is subject to
\begin{align}
\text{s.t.}~&\textstyle\sum_{k \in K} v_{ij}^k \leq y_{ij} M, &\forall i,j, \label{numberassigned}\\
&d_{ij} y_{ij} \leq d_{max}, &\forall i,j, \label{distmax}\\
&\textstyle\sum_{j \in J} v_{ij}^k = w_i^k, &\forall i,k, \label{demandfulfilled}\\
&\textstyle\sum_{i \in I}\sum_{k \in K} t_k v_{ij}^k \leq c_j u_j, &\forall j, \label{charging_cap}\\
&u_j \leq x_j q_j, &\forall j, \label{capacityopened}\\
&\textstyle\sum_{i \in I} y_{ij} \leq x_j M, &\forall j, \label{assignedopened}\\
&\textstyle\sum_{j \in J} y_{ij} \geq 1, &\forall i, \label{demandassigned}\\
&\textstyle\sum_{j \in J} x_j \leq N, & \label{maxstations}\\
&\textstyle\sum_{j \in J}\sum_{l=1}^n \tfrac{1}{n} y_{ij} I_{jl} \geq \bar{p}, &\forall i, \label{reliability}\\
&\textstyle\sum_{j \in J}\sum_{l=1}^n \tfrac{1}{n} v_{ij}^k I_{jl} \geq \bar{p}\textstyle\sum_{j \in J} v_{ij}^k, &\forall i,k. \label{reliability_v}
\end{align}
Constraint \eqref{numberassigned} allows charging only at assigned stations; \eqref{distmax} caps time-to-charge at $d_{max}$; \eqref{demandfulfilled} fulfills all demand ($w_i^k$ vehicles of type $k$ at node $i$); \eqref{charging_cap} respects installed capacity; \eqref{capacityopened} limits connectors per station; \eqref{assignedopened} assigns demand only to opened stations; \eqref{demandassigned} guarantees coverage; \eqref{maxstations} caps the number of stations; and \eqref{reliability}--\eqref{reliability_v} ensure that the probability of at least one assigned station serving a demand without an outage is at least $\bar{p}$, assuming independent outages.

\subsection{Robust vs. Non-Robust Model}
Incorporating $p_j$ boosts robustness against the disruptions that characterize urban areas in developing countries. The Non-Robust Model ignores disruption probability and tends to over-optimize cost, opening too few stations with little coverage redundancy; when an outage hits, demand goes unmet and residents are served poorly. Our formulation is general---setting $I_{jl}=1$ recovers that model---while retaining a linear objective and linearized constraints, so it remains an MIP solvable by standard solvers.

\begin{figure}[t]
    \centering
    \includegraphics[width=\linewidth]{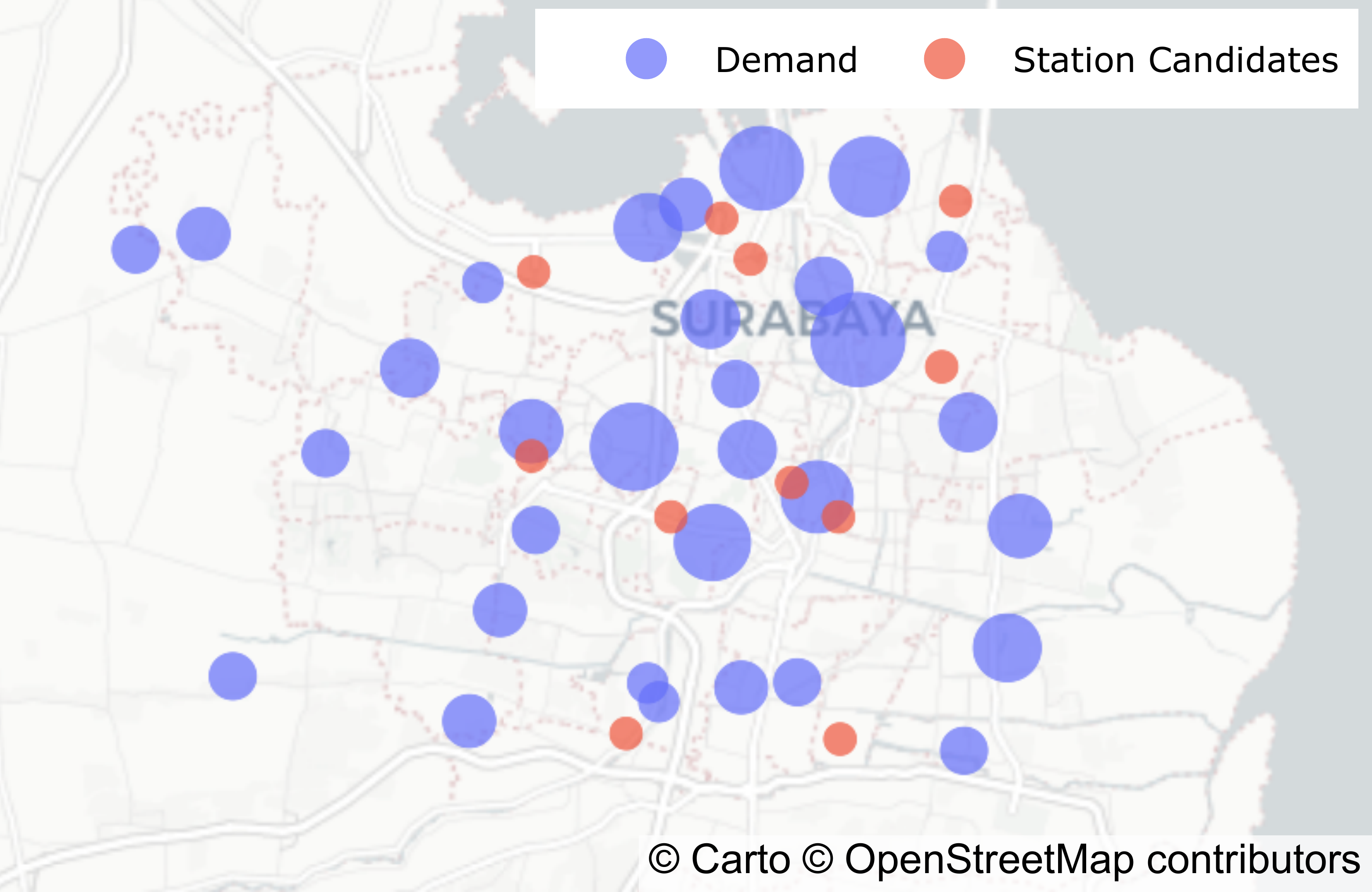}
    \caption{EV charging demand across Surabaya (blue nodes) and the charging-station candidates (red nodes).}
    \label{fig:candidate_location}
\end{figure}

\begin{figure}[t]
    \centering
    \includegraphics[width=0.95\linewidth, trim=0 0.2cm 0 1.1cm, clip]{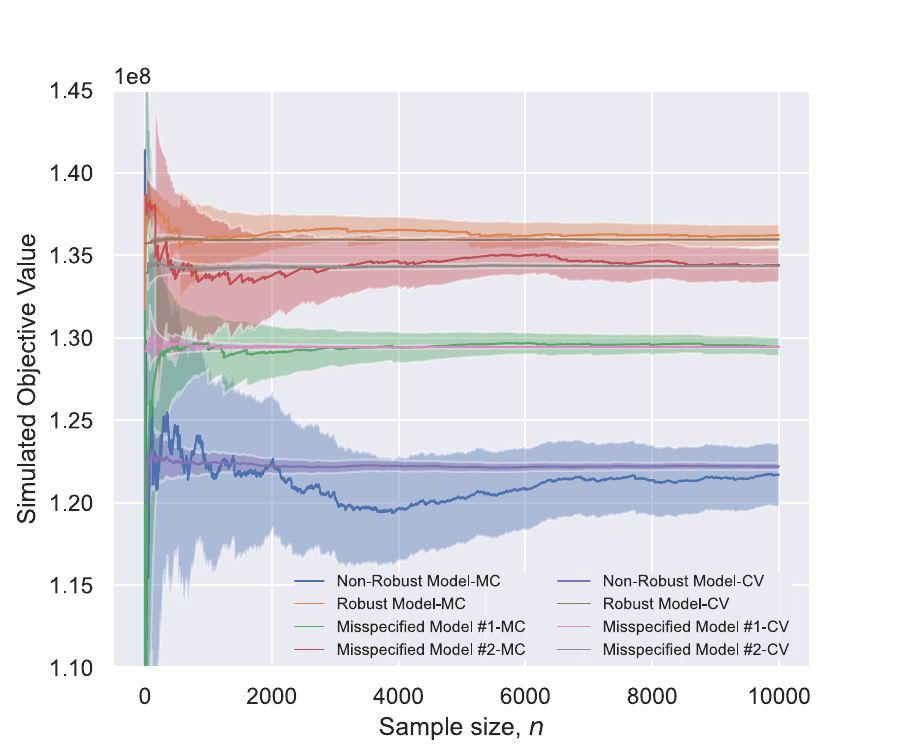}
    \caption{Simulated objective values. Thick lines are means given sample size $n$; shaded areas are 95\% confidence intervals.}
    \label{fig:cv_est}
\end{figure}

\begin{figure*}[t]
    \centering
    \includegraphics[width=0.49\linewidth]{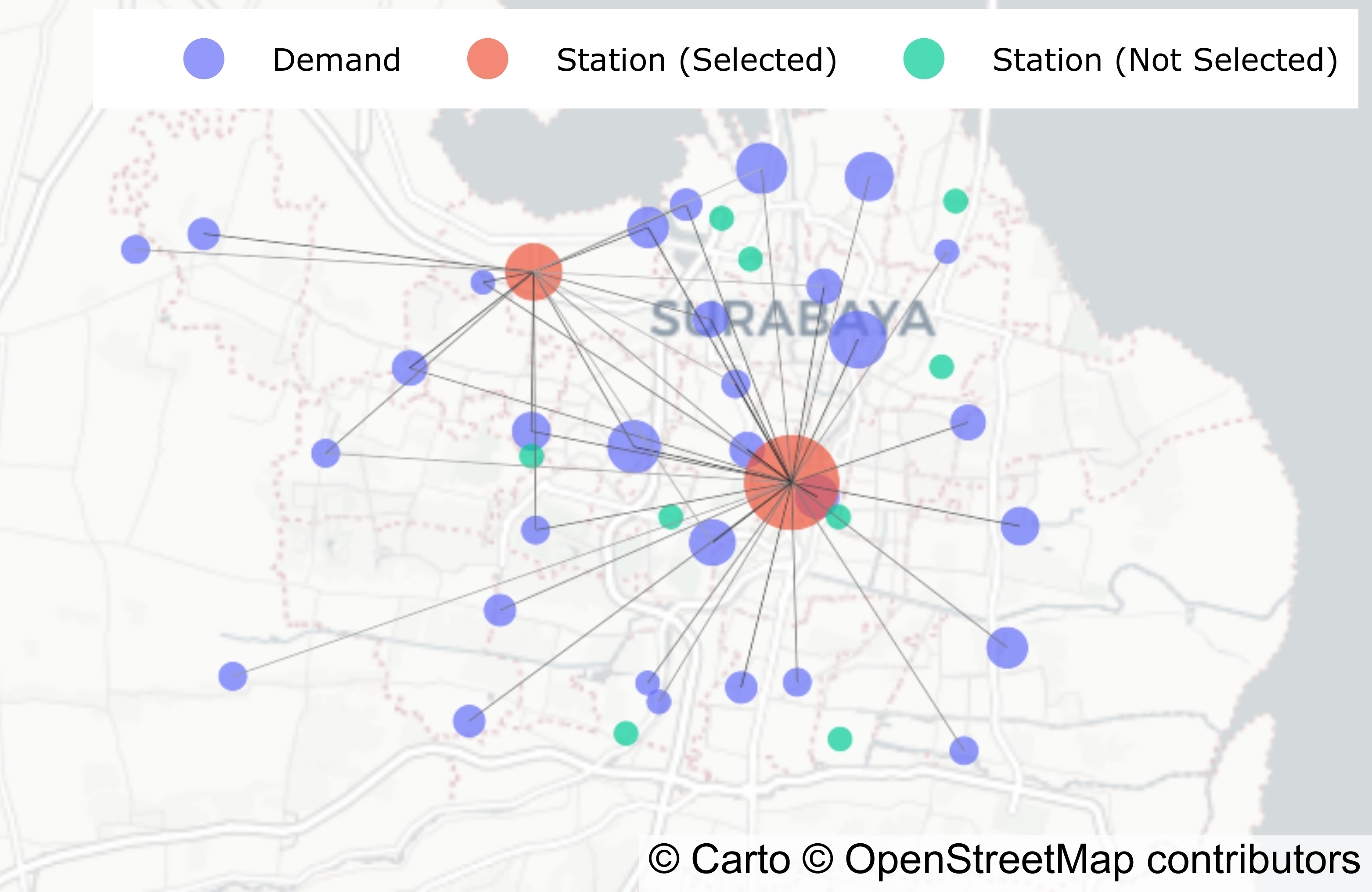}
    \includegraphics[width=0.49\linewidth]{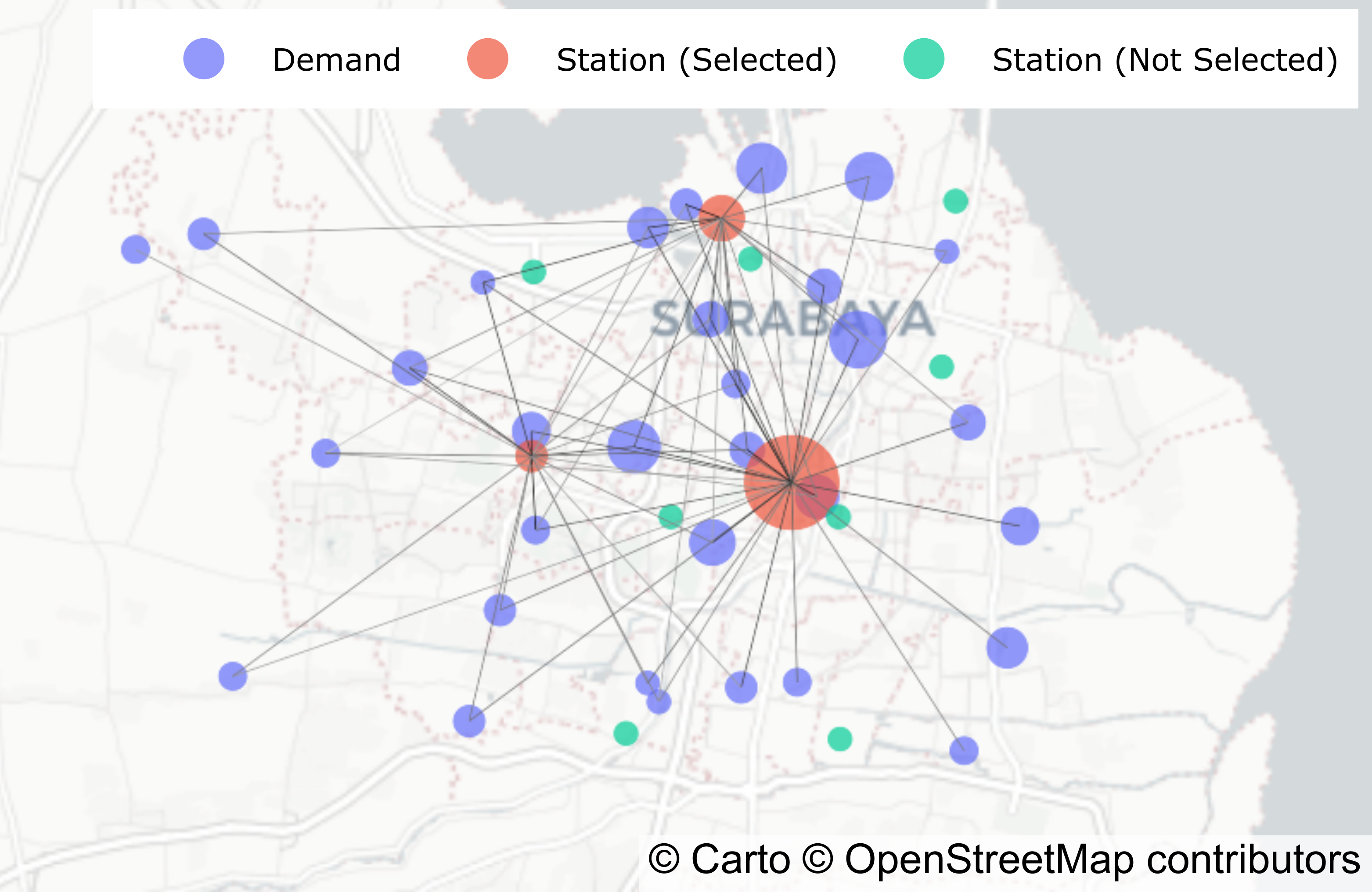}
    \caption{Solutions of the Non-Robust Model (left) and Robust Model (right). The Robust Model distributes capacity to add redundancy against outages.}
    \label{fig:nrm_rm}
\end{figure*}

\subsection{Efficient Reliability Estimation via Control Variates}
The MC estimator $\hat{p}_j = \frac{1}{n}\sum_{l=1}^n \mathbb{I}(Z_{jl} \leq z_j)$ is unbiased and consistent for $p_j$. However, the sample size needed for a precise estimate can be very large when disruptions are rare, motivating variance-reduction and rare-event simulation techniques for reliable certification \cite{arief2021deep}. By Hoeffding's inequality, achieving accuracy $\delta$ at confidence $1-\alpha$ requires $n = \frac{1}{2\delta^2}\ln(2/\alpha)$; e.g., $\delta=10^{-4}$ and $1-\alpha=0.95$ demand $n \approx 1.8\times10^8$ samples \emph{per station}, which can be prohibitive.

To reduce this burden we use control variates (CV) \cite{hesterberg1998control}, introducing an easily computed variable correlated with the quantity of interest. For Gaussian $q=N(\mu,\sigma)$ we use the standard-normal CDF $g(z)=\Phi(z)$ as the control. The CV estimator is
\begin{equation}
\hat{p}_j = \frac{1}{n}\sum_{l=1}^n \mathbb{I}(Z_{jl}\leq z_j) + \pi_j\!\left(\mathbb{I}(X_{jl}\leq \bar{z}_j) - g(\bar{z}_j)\right),
\label{cv_estimator}
\end{equation}
where $X_{jl}$ are standard-normal variables correlated with $Z_{jl}$, $\bar{z}_j$ is the scaled threshold, and $\pi_j$ minimizes the variance:
\begin{equation}
\pi_j = -\frac{\text{Cov}\!\left(\sum_l \mathbb{I}(Z_{jl}\leq z_j),\ \sum_l \mathbb{I}(X_{jl}\leq \bar{z}_j)\right)}{\text{Var}\!\left(\sum_l \mathbb{I}(X_{jl}\leq \bar{z}_j)\right)}.
\end{equation}
The CV estimator \eqref{cv_estimator} is unbiased, and with the optimal $\pi_j$ its variance is reduced by $\text{Cov}^2(\cdot)/[n^2\,\text{Var}(\cdot)] \geq 0$ relative to MC, the reduction being largest when $\mathbb{I}(Z_{jl}\leq z_j)$ and $\mathbb{I}(X_{jl}\leq \bar{z}_j)$ are strongly correlated. This variance reduction lowers the sample size needed for a given $(\delta,\alpha)$, enabling the reliability estimates to be embedded in \eqref{obj_est_sim} and \eqref{reliability}--\eqref{reliability_v} without a large computational penalty.

\section{Numerical Experiments}
\label{sec:exp}

We use EV and electricity data from Surabaya, Indonesia: 11 candidate stations, 31 sub-regions as demand nodes, and two vehicle types (motorcycles $k{=}1$, cars $k{=}2$). Fig.~\ref{fig:candidate_location} shows candidate stations (red) and demand nodes (blue, sized by demand $w_i^k$). Average travel times $d_{ij}$ are from Google Maps. Each connector has $c_j=1440$ min/day (24/7) and stations are limited to $q_j=8$ connectors. We estimate reliability by simulating electricity demand $Z_j \sim q_j=N(\mu_j,\sigma_j^2)$; an outage is triggered when $Z_j \geq z_j$. The masked demand parameters (rescaled by the local utility for privacy) are in Table~\ref{tab:masked_demand_model}, and the remaining parameters in Table~\ref{tab:model_params}.

\begin{table}[t]
\centering
\caption{Electricity Demand $q_j=N(\mu_j,\sigma_j^2)$ and Threshold $z_j$}
\begin{tabular}{cccc}
\toprule
Station $j$ & Mean $\mu_j$ & Std.\ dev.\ $\sigma_j$ & Threshold $z_j$ \\
\midrule
1 & 28{,}979 & 3{,}622 & 36{,}224 \\
2 & 11{,}590 & 1{,}440 & 14{,}490 \\
3 & 5{,}795 & 722 & 7{,}250 \\
4 & 5{,}789 & 715 & 7{,}245 \\
5 & 17{,}387 & 2{,}180 & 21{,}750 \\
6 & 17{,}385 & 2{,}169 & 21{,}730 \\
7 & 11{,}590 & 1{,}442 & 14{,}475 \\
8 & 11{,}595 & 1{,}435 & 14{,}481 \\
9 & 5{,}797 & 720 & 7{,}246 \\
10 & 5{,}785 & 731 & 7{,}241 \\
11 & 11{,}600 & 1{,}450 & 14{,}500 \\
\bottomrule
\end{tabular}
\label{tab:masked_demand_model}
\end{table}

\begin{table}[t]
\centering
\caption{Additional Model Parameters}
\begin{tabular}{llr}
\toprule
Parameter & Symbol & Value \\
\midrule
Station op.+inv.\ cost (range), IDR & $h_j$ & 1.07M--4.60M \\
Electricity price, IDR/kWh & $r$ & 2{,}467 \\
Penalty, IDR/min/vehicle & $s$ & 50{,}000 \\
Charge energy (motorcycle/car), kWh & $e_k$ & 90 / 133 \\
Charge time (motorcycle/car), min & $t_k$ & 20 / 39 \\
Connector op.\ cost, IDR/unit & $g$ & 479{,}285 \\
Max stations & $N$ & 11 \\
Max travel-time threshold, min & $d_{max}$ & 35 \\
Min service-level threshold & $\bar{p}$ & 0.95 \\
\bottomrule
\end{tabular}
\label{tab:model_params}
\end{table}

We build the model by running $n$ simulation replications (up to 10,000) and averaging the objective. We compare four models: the \textit{Robust Model}, the \textit{Non-Robust Model}, and two \textit{Misspecified Models} that assume $0.95p_j$ (underestimation, \#1) and $1.05p_j$ (overestimation, \#2). The CV estimator uses standard-normal $X_{jl}$ with scaled $\bar{z}_j$, giving indicators strongly correlated with $\mathbb{I}(Z_{jl}\leq z_j)$. Fig.~\ref{fig:cv_est} summarizes the simulated objective values with 95\% confidence bands. The code is available at \url{https://github.com/ai-vnv/rare-event-ev-charging}.

\section{Discussion and Findings}
\label{sec:discussion}

\subsection{Robustness of the Optimal Solutions}
Fig.~\ref{fig:cv_est} shows that the Robust Model attains the highest objective, while the Non-Robust Model attains the lowest. By prioritizing cost minimization, the Non-Robust Model opens only two stations, producing low revenue and high penalties during disruptions. The Robust Model instead opens three stations, distributing large stations across the city; this raises total cost by about 18\% but protects against revenue loss and penalties. The net effect is roughly 10\% higher revenue and 60\% lower penalty under disruptions, i.e., about a 13\% higher overall objective. Fig.~\ref{fig:nrm_rm} illustrates that the Robust Model covers more demand points with greater redundancy.

The Misspecified Models retain partial robustness, with objectives between the Robust and Non-Robust solutions: underestimating reliability yields a more conservative four-station solution (over 90\% of demand points covered by two or more stations), while overestimating yields three stations close to the Robust Model. Accurate reliability estimates thus help, but the model tolerates moderate misspecification.

\subsection{Improved Simulation Efficiency via Control Variates}
The central computational challenge is estimating station reliabilities $p_j$ precisely enough to evaluate \eqref{obj_est_sim} and the reliability constraints. While both MC and CV estimators are unbiased and converge to the same values, CV substantially reduces variance, giving the tighter confidence bands in Fig.~\ref{fig:cv_est}. Across estimators, CV attains roughly the same precision with up to $10\times$ fewer samples and no loss of accuracy. This efficiency is what makes embedding the estimated reliabilities in the MIP practical, and it extends the method to settings where disruptions are extremely rare and must be estimated to very high precision (as in developed countries).

\subsection{Applicability and Limitations}
While the results are appealing, our CV estimator  assumes some prior knowledge of the disruption-triggering distribution. Machine-learning models could instead extract disruption features from historical data and estimate EV battery capacity \cite{zhao2023machine} for nonlinear extensions in future work. Finally, in developed countries such as the USA or UK, rare but costly outages make CV-based estimation especially valuable and richer outage records ease calibration. 

\section{Conclusion}
\label{sec:conclusion}

Public EV charging networks in developing countries are being built on power grids where outages still interrupt service, while most siting models assume that every station is always available and therefore yield designs with little redundancy. To address this gap, we propose a simulation-based optimization model for robust EV charging-station planning in developing countries, where frequent electricity disruptions threaten service reliability. The model treats station reliability as a first-class decision factor and uses the control variates variance-reduction technique to estimate it efficiently. On a dataset from Surabaya, Indonesia, the robust solution outperforms its non-robust counterpart by about 13\% and remains robust under reliability misspecification, while CV cuts the required sample size by up to $10\times$ relative to Monte Carlo. The approach offers engineers a practical, disruption-aware decision tool for charging-infrastructure planning.

\bibliographystyle{ieeetr}
\bibliography{ref}

\end{document}